\RequirePackage{fix-cm}
\documentclass[smallextended]{svjour3}     
\smartqed  
\usepackage{graphicx}

\usepackage{amssymb}

\usepackage{pstricks}
\usepackage{pst-node}
\usepackage{url}
\usepackage{subfigure}

\usepackage{dcpic,pictexwd}

\usepackage[authoryear]{natbib}
\usepackage{url}

\begin{document}

\title{Towards a unified framework for decomposability of processes}


\author{Valtteri Lahtinen \and Antti Stenvall
}


\institute{Electromagnetics, Department of Electrical Engineering, Tampere University of Technology, PO Box 692, 33101 Tampere, Finland \\
           Tel.: +358-40-8490430\\
					\url{http://notjargon.org} \\
            \email{valtteri.lahtinen@tut.fi}
						}

\date{Dated: 12.9.2016}

\maketitle

\newcommand{\eqref}[1]{(\ref{#1})}
\newcommand{\fref}[1]{figure~\ref{#1}}
\newcommand{\rmd}{{\rm d}}
\newcommand{\TM}{{\mathrm{T}\Omega}}
\newcommand{\TstarM}{{\mathrm{T}^{*}\Omega}} 

\newcommand{\C}[1]{{\color{red} #1}}
\newcommand{\Cs}[1]{{\color{red}\sout{#1}}}
\newcommand{\Cc}[1]{{\color{red}(Comment: #1)}}

\newtheorem{mydef}{Definition}
\newtheorem{mytheorem}{Theorem}

\begin{abstract}

The concept of process is ubiquitous in science, engineering and everyday life. Category theory, and monoidal categories in particular, provide an abstract framework for modelling processes of many kinds. In this paper, we concentrate on sequential and parallel decomposability of processes in the framework of monoidal categories: We will give a precise definition, what it means for processes to be decomposable. Moreover, through examples, we argue that viewing parallel processes as coupled in this framework can be seen as a category mistake or a misinterpretation. We highlight the suitability of category theory for a structuralistic interpretation of mathematical modelling and argue that for appliers of mathematics, such as engineers, there is a pragmatic advantage from this.

\end{abstract}

\section{Introduction}

The role of category theory as the foundational language for mathematics, or even as \emph{the foundation for mathematics}, has been under discussion for a long time (see e.g. \citep{Lawvere}, \citep{Marquis}, \citep{Muller}). In particular, category theory has been linked to, and seen to allow, a structuralistic interpretation of mathematics, as discussed by \citeauthor{Landry} (\citeyear{Landry}) and \citeauthor{Pedroso} (\citeyear{Pedroso}) among others. Indeed, category theory seems to be more concerned with relations between objects than objects themselves. Moreover, \citeauthor{Arbib} have argued, that not only pure mathematics, but also problems of \emph{applied mathematics}, may be more rigidly founded on a category theoretical soil (\citeyear{Arbib}). In spite of this, the role of category theory as a foundational language in applied mathematics and mathematical modelling is still a rather unexplored path. 

To explore the suitability of category theory for applied mathematics, in this paper we will examine an abstract formal framework for modelling \emph{processes}, utilizing the tools of category theory. Moreover, through our framework and its concrete examples, we argue that the structuralistic view of mathematical modelling implied by our category theoretical approach yields practical benefits for appliers of mathematics. In particular, we will find that \emph{monoidal categories} are a powerful abstraction for the description of parallel and sequential processes. Or put another way, we will here only consider processes that fit into the framework of monoidal categories.\footnote{We note that the abstraction level of monoidal categories might not be adequate for all imaginable processes, and for example higher category theory might be needed.} But as we shall see, this will yield us a fitting framework for some important concrete examples.

\subsection{Motivation}

In science and engineering we rely on modelling the world as more or less interrelated processes. There are \emph{systems} which have \emph{states} and which undergo processes which alter them somehow. Processes are ubiquitous: Be it a chemical process, or the process of aging breaking down to the biological processes happening in our cells, or a computer program building up from different computational processes, that is simply how we perceive the vastness around us. If we did not mentally divide things like this, we would be at a loss with the plethora of information we have to deal with. 

We speak about processes quite freely, and we all probably have an intuitive picture about them. So on an abstract intuitive level, we have an idea about what a process is. Intuitively, a process takes a system from a state to another. It is thus quite natural to think of processes in terms of some kind of diagrams, pictorial representations of them. For example, we could depict a process $f$ taking $A$ to $B$ as

\[\begindc{\commdiag}[50]
\obj(0,10)[A]{$A$}
\obj(0,0)[B]{$B$}
\mor{A}{B}{$f$}
\enddc\]
\noindent or any other intuitive way. This diagram captures the idea of existence of a process between $A$ and $B$: There is a way of getting from $A$ to $B$, described by the arrow $f$. As one may notice, the diagram above is, suggestively enough, how one depicts morphisms in category theory \citep{Adamek}.

In this paper, we are especially interested in \emph{decomposability} of processes. That is, we would like to have a formal framework, which answers the following intuitive question: Can we find two parallel (non-trivial) processes

\[\begindc{\commdiag}[50]
\obj(0,10)[C_1]{$C_2$}
\obj(0,0)[D_1]{$D_2$}
\obj(-10,10)[C_2]{$C_1$}
\obj(-10,0)[D_2]{$D_1$}
\mor{C_1}{D_1}{$g_2$}
\mor{C_2}{D_2}{$g_1$}
\enddc\]

\noindent that put together describe the process

\[\begindc{\commdiag}[50]
\obj(0,10)[A]{$A$}
\obj(0,0)[B]{$B$}
\mor{A}{B}{$f$}
\enddc\]

\noindent exactly? Or, in a similar manner, it may make sense to ask whether the above process $f$ could be described as a sequence
\[\begindc{\commdiag}[50]
\obj(0,10)[A]{$A$}
\obj(0,0)[C]{$C$}
\obj(0,-10)[B]{$B$}
\mor{A}{C}{$f_1$}
\mor{C}{B}{$f_2$}
\enddc\]
\noindent of two consecutive processes $f_1$ and $f_2$. 

But why do we find this interesting? We write from the point of view of engineers, who have found, by trial and error, that to prevent simulation software from quickly becoming obsolete, one must first search for the most profound mathematical foundations of the things one aims to model, when developing such software in order to make it general. Developing software for particular cases only is, of course, often necessary and can lead to very specifically tailored tools, but a more general program can answer to a broader class of needs and perhaps even to needs that no-one has at the present time but will have later on. This has led us to study the formal meaning of processes. For example, in physics and engineering design, the process of solving \emph{field problems} is essential. It is puzzling that so-called \emph{coupled} field problems, or multiphysics problems, are so common in engineering and natural sciences (see \citep{Bossavit}, \citep{Felippa} and \citep{Reddy} for some examples), yet the formal meaning of the term is somewhat unclear. From a structural realist's point of view, it seems that such confusion is a result of not fully recognizing the mathematical structures used in modelling, or utilizing a model not fundamental enough.\footnote{By structural realism, we simply mean realism about structure, as discussed in \citep{North}, implying that formulations of theories are not equivalent if they utilize different structures. For example, North argues that Hamiltonian mechanics is a more fundamental description of reality than Lagrangian mechanics, as it gets by with less structure, even though they result in exactly the same \emph{predictions}.} Irrelevant structure in models obscures their interpretation.\footnote{Consider, for example, the Faraday's law of electromagnetic induction, which states that a time-varying magnetic field induces an electric field, and which can be expressed as $\rmd E = -\partial_t B$, where $E$ is a differential 1-form called the electric field intensity, $B$ is a differential 2-form called the magnetic flux density, $\rmd$ is the exterior derivative and $\partial_t$ the time-derivative operator. This equation requires a differentiable structure from the space but presumes nothing from its metric properties. However, plugging in the metric tensor and representing the same law utilizing the vectorial counterparts of $B$ and $E$ and the metric-dependent vector differential operator $\mathrm{curl}$, as $\mathrm{curl} ({\bf E}) = -\partial_t {\bf B}$, completely overshadows its metric-independent nature, leaving room for misinterpretations: As $\mathrm{curl}$ depends on the metric tensor of the space, it would now seem that electromagnetic induction somehow couples to the metric properties of the space, which, according to Maxwell's theory, it does not. For a physics-oriented introduction to the mathematics of these issues, we refer the reader to \citep{Frankel}.} This is why we want to pursue a reductionistic approach, working towards a general and formal framework for processes, stripped of inessential structure and leaving the instantiations for the modeller.

Coupling is obviously related to decomposability: Can we decompose coupled things to uncoupled ones? Or are coupled things non-decomposable? Given a process, an interesting question thus is, whether there exists a formal framework which would give the tools to separate it into two or more, perhaps simpler ones, if possible. However, an essential part of science is to unify and generalize things, reducing particular concepts to general rules and frameworks. So in the spirit of scientific reductionism, to study the decomposability of processes, we search for a unified abstract framework for the decomposability of $\emph{any}$ process, the concretizations of which would give particular examples of processes. Obviously, this end-goal is not very modest, and this paper will only aim to be the first step towards it.

This work parallels such recent contributions as \citep{CategoriesInControl}, \citep{PassiveLinear}, \citep{Markov}, \citep{RosettaStone} and \citep{Coecke}, all of which are leaps towards \emph{a general science of systems and processes}, and through which category theory has found its way closer to everyday lives of physicists and engineers. However, our work differs from them in the sense that we look at the situation more from a reductionistic engineer's point of view, and we want to discuss and formalize when we can and what it means to \emph{decompose} processes. Moreover, we are especially interested in the conception of parallel processes that are coupled. We argue, that such coupling need not exist when working in the right category.

With this paper, we thus hope to appeal to a two-fold audience. On one hand, we aim to provide insight to mathematically and philosophically oriented engineers and physicists about the nature of decomposable processes through a category theoretical description. On the other hand, we hope that also mathematicians and philosophers with an orientation to applications will find this framework and its concrete examples, as well as their implications, interesting. We emphasize, that this is not a paper written from a philosopher's point of view. We write as engineers, for engineers, physicists, mathematicians and philosophers, to inspire discussion over disciplinary boundaries about the applicability of category theory for modelling processes we confront and the fundamental nature of processes in general.

\subsection{Structure of the article}

We expect the reader to be willing to learn some elementary category theory. In particular, objects, morphisms between objects, functors between categories and natural transformations between functors will be introduced here, and they are essential for our framework. However, no \emph{prior} encounter with category theory is necessarily required. Moreover, we will introduce an extra structure on a category: a monoidal category. However, to keep this work as accessible to general audience as possible, we will not try to be complete in our exposition of monoidal categories in the sense that we would go on to define any extra terminology related to them.\footnote{For example, a monoidal category can be \emph{symmetric}, \emph{closed} or \emph{braided}.} For a good, thorough introduction to category theory, we refer the reader to \citep{Adamek}, and for a more physics-oriented introduction, covering also monoidal categories, we recommend \citep{Coecke}.

In section~2, we give a brief introduction to the basic concepts of category theory, and in section~3, we define what is meant with a monoidal category. Furthermore, in the same section, we describe what it means for processes to be decomposable in this framework. In section~4, we will move on to concretizations of our framework, yielding examples of decomposability of familiar processes. Finally, in section~5, we summarize, draw conclusions and set up the stage for future research in this area.

\section{Category theoretical prerequisites}

In this section, we will define a few fundamental category theoretical concepts. A reader familiar with basic category theory can skip this section.

A \emph{category} may be defined in several equivalent ways. Our definition here is the one presented in~\citep{Adamek}. Intuitively, a category is a bunch of objects with a bunch of arrows between them, and the arrows can be composed to yield new arrows.

\begin{mydef}\label{def:category}
A {\bf category} ${\bf C}$ is a 4-tuple $(\mathcal{O}, \mathrm{hom}_{{\bf C}}, id, \circ)$, which consists of
\begin{itemize}
\item $\mathcal{O}$, a class of {\bf objects},
\item for each pair of objects $(A,B)$ a set $\mathrm{hom}_{{\bf C}}(A,B)$ of {\bf morphisms} from $A$ to $B$,
\item The {\bf identity morphism} $id_A \in \mathrm{hom}_{{\bf C}}(A,A)$ for each object $A$,
\item {\bf composition} $\circ$, which for each $f \in \mathrm{hom}_{{\bf C}}(A,B)$ and for each $g \in \mathrm{hom}_{{\bf C}}(B,C)$ associates a morphism $g \circ f \in \mathrm{hom}_{{\bf C}}(A,C)$,\footnote{The composition of morphisms $g$ and $f$, $g \circ f$, can be conveniently read as \emph{"$g$ after $f$"}.} 
\end{itemize}
such that
\begin{itemize}
\item the composition $\circ$ is associative,
\item for each $f \in \mathrm{hom}_{{\bf C}}(A,B)$, $id_B \circ f = f \circ id_A = f$,
\item The sets of morphisms  $\mathrm{hom}_{{\bf C}}(A,B)$ are pairwise disjoint.
\end{itemize}

\end{mydef}

\noindent Often, the existence of a morphism $f \in \mathrm{hom}_{{\bf C}}(A,B)$ is depicted as an arrow going from $A$ to $B$, similarly as we drew processes in the introduction, or simply writing $f: A \rightarrow B$. We denote the class of morphisms of a category ${\bf C}$ as $\mathrm{Mor}({\bf C})$ and its class of objects as $\mathrm{Obj}$({\bf C}). Objects $A$ and $B$ are said to be {\bf isomorphic} when there exist morphisms $f: A \rightarrow B$ and $g: B \rightarrow A$ such that $g \circ f = id_A$ and $f \circ g = id_B$. The morphism $f$ (and $g$) is then an {\bf isomorphism}. The existence of an isomorphism between objects $A$ and $B$ will be denoted as $A \sim B$. Furthermore, by a {\bf non-trivial morphism}, we will here mean a morphism that is not an identity morphism.

Things get even more interesting, when we start to talk about arrows between categories. These are called \emph{functors}.

\begin{mydef}\label{def:functor}
A {\bf functor} $F$ from a category ${\bf C_1}$ to a category ${\bf C_2}$, denoted $F: {\bf C_1} \rightarrow {\bf C_2}$, is a mapping associating
\begin{itemize}
\item to each object $A \in \mathrm{Obj}({\bf C_1})$ an object $F(A) \in \mathrm{Obj}({\bf C_2})$,
\item to each morphism $f \in \mathrm{hom}_{{\bf C_1}}(A,B)$ a morphism $F(f) \in \mathrm{hom}_{{\bf C_2}}(F(A),F(B))$,
\end{itemize}
such that $F$ preserves identities and composition, i.e.,
\begin{itemize}
\item for each $A \in \mathrm{Obj}({\bf C_1})$, $F(id_A) = id_{F(A)}$,
\item whenever $g \circ f$ is defined, $F(g \circ f) = F(g) \circ F(f)$.
\end{itemize}
\end{mydef}

\noindent Also functors can be composed: Given two functors $F_1: {\bf C_1} \rightarrow {\bf C_2}$ and $F_2: {\bf C_2} \rightarrow {\bf C_3}$ their {\bf composite functor} $F_2 \circ F_1: {\bf C_1} \rightarrow {\bf C_3}$ maps each $A \in \mathrm{Obj}({\bf C_1})$ to $F_2(F_1(A)) \in \mathrm{Obj}({\bf C_3})$ and each $f \in \mathrm{hom}_{{\bf C_1}}(A,B)$ to $F_2(F_1(f)) \in \mathrm{hom}_{{\bf C_3}}(F_2(F_1(A)),F_2(F_1(B)))$. Moreover, the {\bf identity functor} $Id_{\bf C}: {\bf C} \rightarrow {\bf C}$ maps each object and morphism in ${\bf C}$ to itself, and acts thus as an identity with respect to functor composition.

So now we have arrows between objects and arrows between categories, which we call functors. But why stop now? We can further define arrows between functors. These arrows are called \emph{natural transformations}.

\begin{mydef}\label{def:naturalTransformation}
Given two functors $F,G: {\bf C_1} \rightarrow {\bf C_2}$, a {\bf natural transformation} $\eta$ from $F$ to $G$ is a mapping assigning to each $A \in \mathrm{Obj}({\bf C_1})$ a morphism $\eta_A: F(A) \rightarrow G(A)$ such that for each $f \in \mathrm{hom}_{{\bf C_1}}(A,B)$ the diagram
\[\begindc{\commdiag}[50]
\obj(0,10)[a]{$F(A)$}
\obj(10,10)[c]{$G(A)$}
\obj(10,0)[d]{$G(B)$}
\obj(0,0)[b]{$F(B)$}
\mor{a}{c}{$\eta_A$}[\atleft, \solidarrow]
\mor{c}{d}{$G(f)$}[\atleft, \solidarrow]
\mor{a}{b}{$F(f)$}[\atright, \solidarrow]
\mor{b}{d}{$\eta_B$}[\atleft, \solidarrow]
\enddc\]
\noindent commutes\footnote{The \emph{commutation} of the diagram means that $G(f) \circ \eta_A = \eta_B \circ F(f)$.}. The morphisms $\eta_A$ are the {\bf components} of $\eta$.  
\end{mydef}
\noindent A {\bf natural isomorphism} is a natural transformation whose components are isomorphisms.

\section{Monoidal categories: an abstract framework for processes}

Presupposing the definition of a category, in this section we define \emph{monoidal categories}, and some terminology related to them. Then, we go on to discuss what it means to decompose processes in such a setting. 

\subsection{Monoidal categories -- a definition}

In the following, note that the {\bf cartesian product of categories} ${\bf C_1}$ and ${\bf C_2}$, denoted as ${\bf C_1} \times {\bf C_2}$, is the category in which objects are pairs of objects with one element of the pair from each category and morphisms are pairs of morphisms in a similar manner. Identity morphisms and composition are defined for the pairs componentwise.

\begin{mydef}

A {\bf monoidal category} is the 5-tuple $({\bf C}, \otimes, I, \alpha, \lambda, \rho)$, where

\begin{itemize}

\item ${\bf C}$ is a category, 

\item $\otimes: {\bf C} \times {\bf C} \rightarrow {\bf C}$ is the {\bf monoidal product} functor,

\item $I \in \mathrm{Obj} ({\bf C})$ is the {\bf unit object},

\item $\alpha$ is a natural isomorphism called the {\bf associator}, which assigns an isomorphism $\alpha_{A,B,C}:(A \otimes B) \otimes C \sim A \otimes (B \otimes C)$ to each $A,B,C \in \mathrm{Obj} ({\bf C})$,

\item $\lambda$ is a natural isomorphism called the {\bf left unitor}, which assigns an isomorphism $\lambda_A: I \otimes A \sim A$ to each $A \in \mathrm{Obj} ({\bf C})$,

\item $\rho$ is a natural isomorphism called the {\bf right unitor}, which assigns an isomorphism $\rho_A: A \otimes I \sim A$ to each $A \in \mathrm{Obj} ({\bf C})$,

\end{itemize}

\noindent such that the triangle

\[\begindc{\commdiag}[50]
\obj(0,20)[a]{$(A \otimes I) \otimes B$}
\obj(20,20)[b]{$A \otimes (I \otimes B)$}
\obj(20,0)[d]{$A \otimes B$}
\mor{a}{b}{$\alpha_{A,I,Y}$}
\mor{a}{d}{$\rho_A \otimes id_B$}[\atright, \solidarrow]
\mor{b}{d}{$id_A \otimes \lambda_B$}[\atleft, \solidarrow]
\enddc\]

\noindent and the pentagon 
\[\begindc{\commdiag}[50]
\obj(-15,0)[C_1]{$((A \otimes B) \otimes C) \otimes D$}
\obj(15,0)[C_2]{$(A \otimes (B \otimes C)) \otimes D$}
\obj(15,-15)[C_3]{$A \otimes ((B \otimes C) \otimes D)$}
\obj(0,-30)[C_4]{$A \otimes (B \otimes (C \otimes D))$}
\obj(-15,-15)[C_5]{$(A \otimes B) \otimes (C \otimes D)$}
\mor{C_1}{C_2}{$\alpha_{A,B,C} \otimes id_D$}
\mor{C_2}{C_3}{$\alpha_{A,B \otimes C, D}$}
\mor{C_3}{C_4}{$id_A \otimes \alpha_{B,C,D}$}
\mor{C_5}{C_4}{$\alpha_{A,B,C \otimes D}$}[\atright, \solidarrow]
\mor{C_1}{C_5}{$\alpha_{A \otimes B,C,D}$}[\atright, \solidarrow]
\enddc\]

\noindent commute for all $A,B \in \mathrm{Obj} ({\bf C})$, and for all $A,B,C,D \in \mathrm{Obj} ({\bf C})$, respectively.

\end{mydef}

\noindent When it is clear from the context, we may denote a monoidal category \linebreak $({\bf C}, \otimes, I, \alpha, \lambda, \rho)$ simply as ${\bf C}$ or $({\bf C}, \otimes)$. The coherence conditions described by the above diagrams assure the equality of various isomorphisms constructed utilizing the associator and unitors. Remarkably, these conditions are adequate for showing that all isomorphisms, constructed utilizing the unitors and the associator, having the same source and target objects, are the same \citep{MacLane2}.

Note, how the above definition declares e.g. $(A \otimes B) \otimes C$ and $A \otimes (B \otimes C)$ as isomorphic but not necessarily equal, and the same complication applies when operating with unit objects. However, we do not really need to be worried about this, as a monoidal category is always categorically equivalent\footnote{The equivalence of categories guarantees that they behave similarly in terms of categorical properties. Categories ${\bf C_1}$ and ${\bf C_2}$ are {\bf equivalent} provided that there exist functors $F_1: {\bf C_1} \rightarrow {\bf C_2}$ and $F_2: {\bf C_2} \rightarrow {\bf C_1}$ such that $F_2 \circ F_1$ and $F_1 \circ F_2$ are naturally isomorphic to the identity functors $Id_{{\bf C_1}}$ and $Id_{{\bf C_2}}$, respectively. \citep{Adamek}, \citep{Coecke}} to a {\bf strict monoidal category}, in which the associators and unitors are identities \citep{MacLane}. This also allows very descriptive graphical calculi to be developed for monoidal categories \citep{Selinger}.

From a modeller's point of view, monoidality is a \emph{structure} on a category: Choosing how to make a category into a monoidal one can be seen as a modelling decision. Many categories can be made into monoidal categories in several ways. For example, in the category of sets and functions between them, {\bf Set}, the cartesian product and the disjoint union of sets, which are the \emph{product} and \emph{co-product} in {\bf Set}, respectively, are both eligible monoidal products for {\bf Set}. Whenever the monoidal product of a category is chosen to be the category theoretical product or the category theoretical co-product (See e.g. \citep{Adamek} for definitions) and terminal object\footnote{{\bf Terminal object} is such an object, that there exists a unique morphism to it from every object in the category. Dually, {\bf initial object} is such that there exists a unique morphism from it to every other object in the category.} as the unit, the resulting monoidal category is called {\bf cartesian monoidal} or {\bf co-cartesian monoidal}, respectively. From physics' point of view, (co-)cartesian monoidal structures are typically thought to be good for describing classical phenomena, while non-cartesian ones, on the other hand, are suitable for modelling quantum phenomena \citep{RosettaStone}, \citep{Coecke}. It is up to the modeller to decide, what instance of the monoidal structure most appropriately reflects the structure of the part of the world she wants to model.

Note also, how monoidal product on a category is not defined as an internal constituent of the category, but as how two categories are related: \emph{an arrow between categories}. Moreover, the associator and unitors, integral for the monoidal structure, essentially describe how arrows like this are related. This shows how relations take precedence over relata in category theory, highlighting a structuralistic interpretation of mathematical modelling: Models fundamentally describe relations between mathematical objects to give meaning to objects themselves. However, we do not want to take any \emph{ontological} stand here about mathematical objects and structures, but we argue that a pragmatic advantage arises for an applier of mathematics by identifying differences in categories in a structural manner. This will be clarified later on.

\subsection{Decomposability of processes in monoidal categories}

From now on, we shall assume that we are working in a monoidal category. The freedom of instantiating the monoidal structure as we wish equips us with the tools to model a multitude of different processes. In addition, an abstract monoidal category gives us the framework to talk about general abstract processes. In particular, we have an almost ready-made intuitive definition for their decomposability.

In our framework, an abstract {\bf process} is represented by a morphism
\[\begindc{\commdiag}[50]
\obj(0,10)[A]{$A$}
\obj(0,0)[B]{$B$}
\mor{A}{B}{$f$}
\enddc\]
in a monoidal category, taking its object $A$ to $B$. Equivalently, we will sometimes denote such a morphism simply as $f: A \rightarrow B$, as is customary. Morphisms of the form $f: I \rightarrow A$ from the unit object can sometimes be viewed as representing {\bf states} of the abstract system.\footnote{A morphism of the form $f: I \rightarrow A$ can be seen as a representation of a state, as it, in a very formal sense, chooses an element from the object $A$. For example, in the monoidal category of sets and functions, with the cartesian product as the monoidal product, the unit object is the terminal object, a singleton set $1$, to which there is a unique function from every other set. Then, a function $f: 1 \rightarrow A$ maps $1$ to a single element of $A$, and thus in this sense, selects an element of $A$. So even though we often do not explicitly speak about what the objects of certain category contain, it does not mean that we have necessarily lost that information in the process of abstraction: We can recover it using the unit object. \citep{Coecke} As pointed out in \citep{Bain}, this is again an example of how category theory often deals externally with internal constituents of things.} From such a viewpoint, a process represented by a morphism \emph{alters the states} of the system.\footnote{This kind of a view is very useful e.g. in physics where we often deal with the \emph{state space} of a system.} {\bf Sequential processes} are then represented by composable morphisms in this category. Two processes $f_1$ from $A$ to $C$ and $f_2$ from $C$ to $B$ 
\[\begindc{\commdiag}[50]
\obj(0,10)[A]{$A$}
\obj(0,0)[C]{$C$}
\obj(0,-10)[B]{$B$}
\mor{A}{C}{$f_1$}
\mor{C}{B}{$f_2$}
\enddc\]
can thus be described as a single process from $A$ to $B$ through their category theoretical composition $f_2 \circ f_1$. Then, {\bf parallel processes}, which we depict as
\[\begindc{\commdiag}[50]
\obj(0,10)[C_1]{$C_2$}
\obj(0,0)[D_1]{$D_2$}
\obj(-10,10)[C_2]{$C_1$}
\obj(-10,0)[D_2]{$D_1$}
\mor{C_1}{D_1}{$g_2$}
\mor{C_2}{D_2}{$g_1$}
\enddc\]
are represented by the monoidal product of the two morphisms:
\[\begindc{\commdiag}[50]
\obj(0,10)[A]{$C_1 \otimes C_2$}
\obj(0,0)[B]{$D_1 \otimes D_2$}
\mor{A}{B}{$g_1 \otimes g_2$}
\enddc\]

Now, let us define, what it means for these processes to be decomposable. First, consider decomposition into sequential processes:

\begin{mydef}

A process represented by a morphism $f: A \rightarrow B$ in a monoidal category ${\bf C}$ is {\bf sequentially decomposable} if there exist an object $C \in Obj({\bf C})$ and non-trivial morphisms $f_1: A \rightarrow C$, $f_2: C \rightarrow B$, such that the triangle
\[\begindc{\commdiag}[50]
\obj(0,10)[a]{$A$}
\obj(10,10)[c]{$C$}
\obj(10,0)[b]{$B$}
\mor{a}{b}{$f$}[\atright, \solidarrow]
\mor{a}{c}{$f_1$}[\atleft, \solidarrow]
\mor{c}{b}{$f_2$}[\atleft, \solidarrow]
\enddc\]
\noindent commutes.

\end{mydef}

\noindent Sequential decomposability of a process thus means simply, that the morphism in question can be expressed as a composition of two non-trivial morphisms. We deliberately leave the degenerate cases outside the definition, as it does not make much sense to decompose a process to a composition of identity and itself. 

Then, consider parallel decomposition:

\begin{mydef} \label{def:parallel}

A process represented by a morphism $g: C \rightarrow D$ in a monoidal category ${\bf C}$ is {\bf parallel decomposable} if there exist objects $C_1, C_2, D_1, D_2 \in Obj({\bf C})$ such that $C_1 \otimes C_2 \sim C$ and $D_1 \otimes D_2 \sim D$, and non-trivial morphisms $g_1: C_1 \rightarrow D_1$, $g_2: C_2 \rightarrow D_2$ such that the diagram
\[\begindc{\commdiag}[50]
\obj(0,10)[a]{$C_1 \otimes C_2$}
\obj(10,10)[c]{$C$}
\obj(10,0)[d]{$D$}
\obj(0,0)[b]{$D_1 \otimes D_2$}
\mor{a}{c}{$\sim$}[\atleft, 11]
\mor{c}{d}{$g$}[\atleft, \solidarrow]
\mor{a}{b}{$g_1 \otimes g_2$}[\atright, \solidarrow]
\mor{b}{d}{$\sim$}[\atleft, 11]
\enddc\]
\noindent commutes.

\end{mydef}

\noindent The definition states that it is essentially the same thing to go through the process $g_1 \otimes g_2$ as it is to go through $g$, as the initial and final states for these processes are essentially the same.

What motivates this definition? It is the conception of \emph{coupling} of parallel processes. We argue that we do not need to talk about two separate processes that are somehow coupled. That would be a category error. In a suitable monoidal category, any process that is not parallel decomposable can be viewed simply as a single process. That is, when one talks about \emph{coupled processes}, one is actually dealing with a single process, in a certain monoidal category. The key to relieve the confusion is to define the correct category, in which to discuss the processes, by examining the structural differences between categories. On the other hand, this can be seen as an explicit forbiddance of coupled parallel processes. A single non-decomposable process is simply a single process, and two parallel processes are two \emph{separate} parallel processes. If they seem coupled, one is thinking in terms of a wrong category or misinterpreting the processes. 

A monoidal category thus captures the ideas of sequential and parallel processes. It models, in a formal manner, many intuitive aspects of those concepts. The associativity of composition in a category and the associator of a monoidal category allow us to think of any subset of a set of either parallel or sequential processes as a single process, without altering the meaning of the full set of processes. That is, we can parenthesize the processes as we wish. And as we have seen, it allows us to give a formal meaning for the decomposition of processes, as well.

\section{Examples}

As decomposability depends on the instantiation of the monoidal product, it can be seen in a sense as a modelling decision. That is, the modeller has the freedom to model the world as she wishes, and by clever modelling decisions, she might make processes that are non-decomposable in some category, decomposable in another. This is, however, not to say that this would be easy at all, as even categorifying many of our widely used theories and models is highly non-trivial and a subject of contemporary research \citep{CategoriesInControl}, \citep{PassiveLinear}, \citep{Markov}, \citep{RosettaStone}, \citep{Coecke}, \citep{Lal}.

In this section, we will take a look at decomposability of processes in a few particular monoidal categories. First, we shall familiarize ourselves with the definitions by looking at the simple example of sets and functions between them. Moreover, we argue that \emph{parallel coupling} of processes need not happen, given the correct category for modelling them. Finally, we will discuss quantum mechanical processes in the light of our framework.

\subsection{Sets and functions}

The category ${\bf Set}$ of sets and functions between them becomes a (cartesian) monoidal category by instantiating the cartesian product $\times$ of sets as the monoidal product with a singleton as its unit object and the required natural isomorphisms accordingly. In $({\bf Set}, \times)$, it is also easy to see why we need, e.g., the associators in the definition of monoidal category: even though $(A \times B) \times C \neq A \times (B \times C)$, which is reflected elementwise by the fact that $((a,b),c) \neq (a,(b,c))$, the two sets are indeed isomorphic. So even this simple example does not yield a strict monoidal category but the associators are necessary.

Let us now, however, take the disjoint union $\oplus$ of sets as the monoidal product for ${\bf Set}$ and the empty set $\emptyset$ as the unit object. Here, we will consider that the objects of ${\bf Set}$, which are simply sets, contain the states as their elements, and functions between sets, the morphisms of ${\bf Set}$, will thus be our processes. Obviously, sequentially decomposing a process in this category means simply to express a function as a composition of two functions. But what about parallel decomposability? Consider a process represented by a function
\[\begindc{\commdiag}[50]
\obj(0,10)[A]{$A$}
\obj(0,0)[B]{$B$}
\mor{A}{B}{$f$}
\enddc\]
\noindent where $A = \{a_1,a_2\}$, $B = \{b_1,b_2\}$, and $f(a_1) = b_1$ and $f(a_2) = b_2$. Then consider sets $A_1 = \{a_1\}$, $A_2 = \{a_2\}$, $B_1 = \{b_1\}$ and $B_2 = \{b_2\}$, which also yield the sets $A_1 \oplus A_2 = \{(a_1, 1), (a_2, 2))\}$ and $B_1 \oplus B_2 = \{(b_1, 1),(b_2, 2)\}$ through the chosen monoidal product. Now, obviously $A \sim A_1 \oplus A_2$ and $B \sim B_1 \oplus B_2$ by associating $a_i$ and $b_i$ with the pairs containing $a_i$ and $b_i$, respectively. Furthermore, there exist functions $f_1: A_1 \rightarrow B_1$ and $f_2: A_2 \rightarrow B_2$ with $f_1(a_1) = b_1$ and $f_2(a_2) = b_2$. Hence, $f_1 \oplus f_2: A_1 \oplus A_2 \rightarrow B_1 \oplus B_2$ maps $(f_1 \oplus f_2)((a_1, 1)) = (b_1, 1)$ and $(f_1 \oplus f_2)((a_2, 2)) = (b_2, 2)$. Thus, we have a commutative diagram of the form
\[\begindc{\commdiag}[50]
\obj(0,10)[a]{$A_1 \oplus A_2$}
\obj(10,10)[c]{$A$}
\obj(10,0)[d]{$B$}
\obj(0,0)[b]{$B_1 \oplus B_2$}
\mor{a}{c}{$\sim$}[\atleft, 11]
\mor{c}{d}{$f$}[\atleft, \solidarrow]
\mor{a}{b}{$f_1 \oplus f_2$}[\atright, \solidarrow]
\mor{b}{d}{$\sim$}[\atleft, 11]
\enddc\]
\noindent and $f$ is indeed parallel decomposable in $({\bf Set}, \oplus)$.

The above result is very intuitive. It simply shows that instead of a function that maps $a_1$ to $b_1$ and $a_2$ to $b_2$ we can have two separate functions, the first of them mapping $a_1$ to $b_1$ while the second one maps $a_2$ to $b_2$. This is not necessarily yet the best example of modelling parallel decomposability of real world processes, but it does provide some insight to the significance of definition~\ref{def:parallel}.

\subsection{The misconception of parallel coupling}

Processes that seem to be in parallel, can sometimes be interpreted as being \emph{coupled} together in some way. However, we argue that such an interpretation is typically not necessary. \emph{Parallel coupling can be viewed as a category mistake or a misinterpretation.} That is, the modeller making this interpretation has not either properly defined, which category she is working in, or she could view the seemingly coupled processes as a single process in some other category. However, we do not mean, that viewing processes as coupled could not sometimes be beneficial. But it can be beneficial to think of them as a single one, as well. For example, concerning magnetism and elasticity, Bossavit (\citeyear{Bossavit}) argued that coupling of these phenomena occurs when the morphism $\Psi$ representing the energy density of the system is not a sum of distinct energy density terms $\Psi_\mathrm{mag}$ and $\Psi_\mathrm{ela}$ related to magnetism only and elasticity only, respectively. That is, the description of the process, which associates the magnetic state and the elastic state of the system with an energy distribution, should not be viewed anymore as two different morphisms, but as a single non-decomposable one, which associates the \emph{magneto-elastic state} of the system with an energy distribution. The lack of parallel decomposability for this morphism thus allows us to identify that the phenomena cannot be separated, at least in the representational scheme used. Explicitly considering this example in the framework of monoidal categories could also clarify this viewpoint.

Next, we will give some examples.

\subsubsection{Misinterpreting sequential decomposability as parallel coupling}

Consider the functions $f,g: \mathbb{R} \rightarrow \mathbb{R}$
\begin{equation}\label{f}
f(x) = x^2,
\end{equation}
\begin{equation}\label{g}
g(x) = f^{3}(x). 
\end{equation}
Now, it would seem tempting to say that $f$ and $g$ represent parallel processes that are somehow \emph{coupled} through \eqref{g}. However, a closer inspection reveals, that there is no need for that kind of interpretation. In {\bf Set}, they actually form together a single process
\begin{equation}\label{gf}
(g \circ f)(x) = g(f(x)),
\end{equation}
where
\begin{equation}\label{f2}
f(x) = x^2,
\end{equation}
\begin{equation}\label{g2}
g(x) = x^3. 
\end{equation}
That is, they form a single, sequentially decomposable process. This is a simple manifestation of our claim, that identifying the right category for modelling the processes, coupled parallel processes do not exist.

\subsubsection{Parallel coupling as a category mistake}

Consider the pair of equations
\begin{equation}\label{pair1}
x+y = 3,
\end{equation}
\begin{equation}\label{pair2}
y-x = 1,
\end{equation}
forced to hold for some $x,y \in \mathbb{R}$. Now, clearly, the equations are coupled: there is no way of solving \eqref{pair1} without solving \eqref{pair2} at the same time. So, it might seem they would describe coupled parallel processes.

But again, we do not \emph{need} to look at this pair of equations of coupled parallel processes of some kind. The category ${\bf FinVec}$ has finite-dimensional vector spaces as objects and linear mappings between them as morphisms, and it can be made monoidal e.g. by instantiating the direct sum of vector spaces $\oplus$ as the monoidal product. The pair of equations \eqref{pair1}, \eqref{pair2} can be seen as describing the operation of a morphism in ${\bf FinVec}$. That is, it says that there is a single morphism $\mathbb{R}^2 \rightarrow \mathbb{R}^2$ in ${\bf FinVec}$, which takes the state $(x,y)$ to state $(3,1)$ in a fixed basis. Obviously, in this basis, the morphism is described by the matrix \[ \left( \begin{array}{cc} 1 & 1 \\ -1 & 1 \end{array} \right).\] Then, solving for \eqref{pair1} and \eqref{pair2} accounts for finding an inverse for the morphism represented by this matrix. 

So again, viewing \eqref{pair1} and \eqref{pair2} as coupled parallel processes can be seen as thinking in terms of a wrong category, or more appropriately, not clearly defining the category, in which one is working. The abstract framework allows us to reduce the seemingly coupled processes to a single one.

\subsection{Quantum processes}

For a physicist, perhaps the prime example of (the lack of) decomposability comes from quantum mechanics. The concept of \emph{quantum entanglement} is responsible for many of the counter-intuitive effects of the quantum realm, such as quantum teleportation and the dazzling power of quantum computation \citep{Kvanttikirja1}, \citep{Kvanttikirja2}. Let us now consider entanglement and parallel decomposability of quantum processes in the light of our framework.

The monoidal category, in which quantum physics and computation is naturally modelled, is $({\bf Hilb}, \otimes)$, which has Hilbert spaces\footnote{Recall that a Hilbert space is a vector space equipped with an inner product, complete with respect to the norm induced by that inner product.} as its objects, bounded linear operators between Hilbert spaces as its morphisms, and the tensor product of Hilbert spaces $\otimes$ as the monoidal product. The unit object of this monoidal category is a 1-dimensional Hilbert space. In this description of reality, quantum entanglement is then a reflection of the fact that the monoidal product of $({\bf Hilb}, \otimes)$ is non-cartesian \citep{RosettaStone}, \citep{Coecke}. This shows us again the power of structuralistic, category theoretical thinking: It allows us to identify entanglement as a structural property of $({\bf Hilb}, \otimes)$.

In quantum computation, one talks about {\bf qubits}, which are simply quantum mechanical systems with two-dimensional state spaces. It is typical to denote an orthonormal basis for such a state space as $\{|0\rangle, |1\rangle\}$. The elements $|\psi\rangle$ of the state space are further restricted by the condition $\langle \psi | \psi \rangle = 1$, which is a customary way to write that the inner product of $|\psi\rangle$ with itself equals one. So basically, the elements of the Hilbert space $H \in \mathrm{Obj}({\bf Hilb})$ can be considered as describing the same state, whenever they are vectors in the same phase, as only unit vectors are allowed.

As the state space of a qubit is described by a Hilbert space, a compound system of two qubits, a 2-qubit, is described by the tensor product of the state spaces of the individual qubits $H_1 \otimes H_2$. However, not every element of $H_1 \otimes H_2$ is a tensor product of elements of $H_1$ and $H_2$. Such elements are called {\bf entangled states}. For example, expressed using the basis $\{|0\rangle, |1\rangle\}$, $|0\rangle \otimes |1\rangle$ is indeed clearly a tensor product of two 1-qubit states, but the so-called Bell state $(|0\rangle \otimes |0\rangle + |1\rangle \otimes |1\rangle)/\sqrt{2}$ is not. Hence, it is an entangled state.\footnote{The 1-dimensional Hilbert space $\mathbb{C}$ is the unit object of $\bf{Hilb}$. For such an object it holds that $\mathbb{C} \sim \mathbb{C} \otimes \mathbb{C}$. Hence, 2-qubit states in $H_1 \otimes H_2$ can be seen as morphisms $f: \mathbb{C} \rightarrow H_1 \otimes H_2$, obviously with the restriction to unit vectors described above. Thus, we can actually interpret entangled states through our framework of processes: If a morphism $f: \mathbb{C} \rightarrow H_1 \otimes H_2$ representing a 2-qubit state is decomposable to $f_1 \otimes f_2$, with $f_1: \mathbb{C} \rightarrow H_1$ and $f_2: \mathbb{C} \rightarrow H_2$, the state is not entangled. Otherwise, it is entangled and must be represented as a single non-decomposable morphism.}

But how does entanglement of states affect the decomposability of quantum \emph{processes}? Consider a bounded linear operator $h \in \mathrm{Mor}({\bf Hilb})$ of the form $h = f \otimes g$, such that $f: H_1 \rightarrow H_1$ and $g: H_2 \rightarrow H_2$, and thus $h: H_1 \otimes H_2 \rightarrow H_1 \otimes H_2$. Then, by definition~\ref{def:parallel}, it is obvious that $h$ is decomposable into two parallel processes $f$ and $g$, even though $H_1 \otimes H_2$ will, in general, contain entangled states. The operation of $f \otimes g$ on entangled states is simply linearly extended from its operation on unentangled ones. Hence, the entanglement of states does not directly affect the decomposability of a process. However, it is well-known that not all morphisms of type $H_1 \otimes H_2 \rightarrow H_1 \otimes H_2$ in ${\bf Hilb}$ are tensor products. Consider, for example, the operator $s: \mathbb{C}^2 \otimes \mathbb{C}^2 \rightarrow \mathbb{C}^2 \otimes \mathbb{C}^2$, $\mathrm{s}(|\psi_1\rangle \otimes |\psi_2\rangle) = |\psi_2\rangle \otimes |\psi_1\rangle$. Now, we cannot parallel decompose $s$ in $({\bf Hilb}, \otimes)$. The initial and final state spaces are isomorphic with a tensor product of objects, for example through the identity morphism. However, there are no morphisms $f$ and $g$ such that $id_{\mathbb{C}^2 \otimes \mathbb{C}^2} \circ s = (f \otimes g) \circ id_{\mathbb{C}^2 \otimes \mathbb{C}^2}$. Note however, that again, we do not have two coupled processes, but a single process, that is not parallel decomposable.

\section{Conclusions and outlook}

The concept of \emph{process}, ubiquitous in science, engineering and everyday life, obtains a formal meaning as a morphism in a monoidal category. Furthermore, monoidal categories provide a framework for sequential and parallel \emph{decomposability} of processes. This provides insight to the abstract mathematical structures behind representations of systems with parallel and sequential processes. Furthermore, decomposability is a key issue in, for example, parallelization of computation, when designing software for modelling processes. Thus, a formalization of this concept is necessary. 

In this paper, we have characterized decomposability of abstract processes utilizing the framework of monoidal categories. We gave a definition for a monoidal category, defined its morphisms to represent processes, and described what it means to decompose a process in such a framework. Then, we provided some examples related to modelling. The framework proved to be sufficient to describe decomposability of processes, modelling for example quantum phenomena. Our approach manifests the suitability of category theory for problems of applied mathematics on a very general, abstract level, instantiable to a multitude of concrete, particular cases.

We argued that modelling \emph{coupled parallel processes} is an unnecessary complication, at least on an abstract level: It is a matter of choosing a suitable monoidal category, in which the modelling is performed. That is, talking about coupled processes can be ultimately viewed as a category mistake or a misinterpretation. This view implies that it does not necessarily make sense to look at processes that are inseparably intertwined as different ones, but as a single, non-decomposable process. In terms of, for example, engineering science, in which coupled processes are ubiquitous, this reductionistic approach suggests that engineers could clarify their modelling workflow by explicitly considering the most suitable categories to work in. That is, there is a pragmatic advantage to thinking applied mathematics in terms of categories, as long as one is willing to take the structuralistic view to modelling, associated with the category theoretical approach.

This paper provides a starting point for future research and discussion. Monoidal categories have been shown to be a fitting framework for a multitude of processes, having their applications from chemistry to computation and from electrical networks and control theory to quantum physics and logic. Interesting research questions arising from this treatment are numerous: For example, can we decompose certain so far non-decomposable processes by wise modelling decisions? Or can we find monoidal categories that provide insight to decomposing or composing familiar processes? Moreover, which familiar processes that we often view as coupled, are actually a single non-decomposable process in some monoidal category? From an ontological point of view, on what level does the highly general applicability of monoidal categories in modelling the world reflect the structure of the world? Answering even some of these questions could provide invaluable information to a large number of fields of science, in addition to being of interest on their own right.

\section*{Acknowledgment}
This research was supported by The Academy of Finland project [287027]. We would also like express our gratitude to editor and anonymous reviewers for their helpful comments, which led to many improvements in the manuscript.

\end{document}